\RequirePackage[reqno]{amsmath}
\documentclass[a4paper,12pt]{amsart}

\usepackage{amssymb,mathtools}
\usepackage[reqno,fleqn,intlimits]{empheq}
\renewcommand{\div}{\mbox{div}}

\usepackage{cite}

\DeclareFontFamily{OT1}{rsfs}{}
\DeclareFontShape{OT1}{rsfs}{m}{n}{ <-7> rsfs5 <7-10> rsfs7 <10-> rsfs10}{}
\DeclareMathAlphabet{\mathscr}{OT1}{rsfs}{m}{n}

\newcommand{\bel}[1]{\begin{equation}\label{#1}}
\newcommand{\beal}[1]{\begin{eqnarray}\label{#1}}
\newcommand{\beadl}[1]{\begin{deqarr}\label{#1}}
\newcommand{\eeadl}[1]{\arrlabel{#1}\end{deqarr}}
\newcommand{\eeal}[1]{\label{#1}\end{eqnarray}}
\newcommand{\eead}[1]{\end{deqarr}}
\newcommand{\eea}{\end{eqnarray}}
\newcommand{\eeaa}{\end{eqnarray*}}

\newcommand{\be}{\begin{equation}}
\newcommand{\ee}{\end{equation}}

\DeclareFontFamily{OT1}{rsfs}{}
\DeclareFontShape{OT1}{rsfs}{m}{n}{ <-7> rsfs5 <7-10> rsfs7 <10->
rsfs10}{} \DeclareMathAlphabet{\mycal}{OT1}{rsfs}{m}{n}

\newcounter{mnotecount}[section]

\newcommand{\rmnote}[1]{}

\def\mysavedown#1{\edef\mysubs{\mysubs#1}}
\def\mysaveup#1{\edef\mysups{\mysups#1}}
\def\mydown#1{{\mytensor}_{\vphantom{\mysubs}#1}}
\def\myup#1{{\mytensor}^{\vphantom{\mysups}#1}}
\def\tensor#1#2{
  #1
  \def\mytensor{\vphantom{#1}}
  \def\mysubs{\relax}
  \def\mysups{\relax}
  \let\down=\mysavedown
  \let\up=\mysaveup
  #2
  \let\down=\mydown
  \let\up=\myup
  #2
  }

\newcommand{\Tr}{\operatorname{Tr}}

\newcommand{\R}{\mathbb R}

\renewcommand{\div}{\operatorname{div}}

\renewcommand{\epsilon}{\varepsilon}
\renewcommand{\hat}{\widehat}

\def\crn#1#2{{\vcenter{\vbox{
        \hbox{\kern#2pt \vrule width.#2pt height#1pt
           }
          \hrule height.#2pt}}}}

\newcommand{\gtw}{{\widetilde g}}
\newcommand{\phitw}{{\widetilde \phi}}

\renewcommand{\hbar}{{\overline h}}

\newcommand{\pre}[2]{{{\vphantom{#2}}^{#1}}\kern-.2ex{#2}}

\theoremstyle{plain}
\newtheorem{theorem}{Théorème}[section]

\theoremstyle{definition}

\newtheorem{remarks}[theorem]{Remarks}

\numberwithin{equation}{section}

\date{November 21, 2016}

\begin{document}
\title[Conformal parametrizations of modified gravity initial data ]
{ Conformally covariant  parametrizations  for initial data in some modified Einstein gravity theories}
\author[E. Delay]{Erwann Delay}
\address{Erwann Delay,
Laboratoire de Mathématiques d'Avignon (EA 2151), F-84018 Avignon, France}
\email{Erwann.Delay@univ-avignon.fr}
\urladdr{http://www.univ-avignon.fr/fr/recherche/annuaire-chercheurs\newline$\mbox{ }$\hspace{3cm}
/membrestruc/personnel/delay-erwann-1.html}

\begin{abstract}
 We propose further conformal parametrizations  for initial data
in some modified Einstein gravity theories.
Some of them give rise to conformally covariant systems.
\end{abstract}

\maketitle

\\
\newline
\section{Introduction}
 Motivated by the problem of dark matter, early-time inflation or late-time acceleration
of the universe,
a multitude of modified gravity theories are  currently studied: 
$f(R)$, Lovelock, Einstein-Gauss-Bonnet, cubic, quadratic,...

As in  general relativity, the space and time correlation in theses theories  implies that the initial data for the evolution problem cannot be chosen freely.

In general, in  vacuum for instance, the initial data for such modified Einstein gravity theories are given by a manifold $M$ of dimension $n$, equipped with a Riemannian metric $\widehat{g}$ (spatial geometry at a time $0$)
and a symmetric 2-tensor field $\widehat{K}$ (infinitesimal deformation of the spatial geometry  at time $0$), satisfying  some constraint equations  of the form (for details, see eg. 	\cite{Carges}, \cite{BartnikIsenberg2004}, \cite{LeflochMafR}, \cite{ToriiShinkaiEGBgravity}, and  the appendix \ref{sectionequation} )
\begin{equation}\label{eqConstraints}\tag{C}
\left\{
\begin{aligned}
\rho(\widehat{g},\widehat{K})=0\;, \\
J(\widehat{g},\widehat{K}) =0\;.
\end{aligned}
\right.
\end{equation}
In this system,  the first equation, the Hamiltonian constraint, is scalar   and the second one, the momentum constraint, is vectorial.
Such a system is known as the Cauchy problem for the related gravity theory.
This system is highly under-determined because it contains $(n+1)$ equations for $n(n+1)$
unknowns. As a consequence, it is natural to look for $(n+1)$  unknowns, fixing the remaining ones.\\

The goal of this note is to point out some different possibilities to parametrize
(and then different ways to construct)  solutions of \eqref{eqConstraints}. Some of them give rise to 
an interesting mathematical property, namely the   conformal covariance. \\

\noindent{\sc Acknowledgements}: I am grateful to P.T. Chru\'sciel and Ph. Delano\"e   for comments.
\section{The usual conformal method}
On a smooth manifold $M$ of dimension $n$, given a Riemannian metric $g$, we denote by
 $\nabla$ its Levi-Civita connection. If $h$ is a symmetric covariant 2-tensor field, we define its divergence
as the 1-form given by
$$
(\div_g h)_i=-\nabla^kh_{ki}.
$$ 
The classical method (York's method A)
starts from a given metric $g$ together with a trace-free and divergence-free symmetric 2-tensor field $\sigma$
(a TT-tensor) and a real function\footnote{playing morally the role of a mean curvature function} $\tau$. It consists in looking  for solutions of \eqref{eqConstraints} of the form 
\begin{equation}\label{eqParameterization}\tag{P}
\widehat{g}=\phi^{N-2}g\;,\;\; \widehat{K}=\frac{\tau}n\widehat{g}+\phi^{-2}(\sigma+\mathring {\mathcal L}_g W)
\end{equation}
where $N=\frac{2n}{n-2}$, and the unknowns are a function  $\phi>0$ and a one form $W$ and where
$$
(\mathring{\mathcal L}_g W)_{ij} =\nabla_iW_j+\nabla_jW_i-\frac2n\nabla^kW_k \;g_{ij}.
$$
We infer from \eqref{eqConstraints} and \eqref{eqParameterization} a coupled system
of the form (see eg. \cite{BartnikIsenberg2004}, \cite{ToriiShinkaiEGBgravity})
\begin{subequations}\label{eqConformalA}
\begin{empheq}[left=(S)\empheqlbrace]{align}
\rho(\hat g,\hat K)=:L_{g,\tau,\sigma}(\phi,W)=&0 \label{eqLichnerowicz0}\tag{L}\\
J(\hat g,\hat K)\label{eqVector0}=:V_{g,\tau,\sigma}(\phi,W)=&0\tag{V}
\end{empheq}
\end{subequations}
where the scalar equation \eqref{eqLichnerowicz0} is a generalisation  of the Lichnerowicz one (see \cite{Lichnerowicz:44})  and the  equation \eqref{eqVector0} is usually called the  vector equation (see \cite{York:decompo}). \\

In order to produce further interesting parametrization we recall some basic fact. 
Firstly, the operators $\div_g$ and $\mathring {\mathcal L}_g$ are conformally covariant (see appendix \ref{cc}
for a precise definition). Indeed,  for any positive  function $\psi$, if $\gtw=\psi^{N-2}g$ then
\begin{equation}\label{eqConformalVector}\tag{$\mathcal{C}$}
\div_{\gtw} h=\psi^{-N}\div_g(\psi^2 h)\;,\;\;\mathring {\mathcal L}_{\gtw}W=\psi^{N-2}\mathring {\mathcal L}_g(\psi^{2-N}W).
\end{equation}

Secondly, we recall the York decomposition \cite{York:decompo} valid for instance if $M$ is compact and $g$ has no
conformal Killing vector fields (i.e. $\ker\;\mathring{\mathcal L}_g$ is trivial): any covariant symmetric trace free 2-tensor
field $h$ splits in a unique way as 
\begin{equation}\label{eqYork}\tag{Y}
 h=\sigma+\mathring{\mathcal L}_gW,
\end{equation}
where $\sigma$ is a TT-tensor and $W$ is a 1-form.\\

\section{A  conformally covariant parametrization}\label{paramCC}
Getting back to \eqref{eqConstraints}, we use now  the York decomposition 
relative to  $\widehat{g}$, namely
\begin{equation}\label{eqhatYork}\nonumber
\widehat{K}-\frac{\tau}n\widehat{g}=\widehat{\sigma}+\mathring {\mathcal L}_{\widehat{g}}\widehat{W}.
\end{equation}
From \eqref{eqConformalVector}, if we are looking for a solution of \eqref{eqConstraints}  in a conformal class $\widehat{g}=\phi^{N-2}g$,
it is natural to introduce $\hat\sigma=\phi^{-2}\sigma $,
where $\sigma$ is a TT-tensor.  
Still using   \eqref{eqConformalVector}, and expressing $\mathring {\mathcal L}_{\widehat{g}}\hat W$  in terms  of $g=\phi^{2-N} \widehat{g}$,
we are also prompted to set $\hat W\coloneqq\phi^{N-2} W$.

Sticking to the same fixed  $g,\sigma,\tau$, we can thus  parametrize the solutions of the constraint \eqref{eqConstraints} by
\begin{subequations}\label{eqConformalBparam}
\begin{empheq}[left=\empheqlbrace]{align*}
\widehat{g}&=\phi^{N-2}g\;,\\
\widehat{K}&=\frac{\tau}n\widehat{g}+\phi^{-2}(\sigma+\phi^N\mathring {\mathcal L}_g W)\;.
\end{empheq}
\end{subequations}
With this parametrization in \eqref{eqConstraints}
we obtain a new system of the form (see eg. \cite{Delay:paramCC})
\begin{subequations}\label{eqConformalB}
\begin{empheq}[left=(S'_{g,\tau,\sigma})\empheqlbrace]{align*}
\rho(\hat g,\hat K)=:L'_{g,\tau,\sigma}(\phi,W)&=0\;, \\
J(\hat g,\hat K)=:V'_{g,\tau,\sigma}(\phi,W)&=0\;. 
\end{empheq}
\end{subequations}
If we follow the same method but, instead of $g$, starts from a conformal metric $\gtw=\psi^{N-2}g$ and from the related symmetric 2-tensor $\widetilde{\sigma}=\psi^{-2}\sigma$
(which is div$_\gtw$ free by \eqref{eqConformalVector}), sticking to the same given real function $\tau$, we see from \eqref{eqConformalVector} that the couple $(\phitw=\psi^{-1}\phi, \widetilde{W}=\psi^{N-2}W)$ solves the resulting system $(S'_{\widetilde g,\tau,\widetilde\sigma})$ iff $(\phi,W)$ solves $(S'_{g,\tau,\sigma})$.

In other words,the system $(S'_{g,\tau,\sigma})$  is
{\it conformally covariant}.
This simple observation shows that  York's method B (see \cite[Section 4.1]{BartnikIsenberg2004} or the original paper \cite[Page 461]{York:decompo}), also called the physical TT method,
stays valid for every modified gravity theory, regardless of the explicit form of the equations \eqref{eqConstraints}.

\section{Further parametrizations }\label{autreparam}
In this section, we propose  more general but still natural  parametrizations.
Let  $\gtw=\psi^{N-2}g$ be  another conformal metric.
The York  decomposition \eqref{eqYork} relative to $\gtw$ of the tensor
$$\phi^2\psi^{-2}( \widehat{K}-\frac{\tau}n\widehat{g}),$$
leads to the parametrization
$$
\left\{
\begin{aligned}
\widehat{g} &=\phi^{N-2}g\;,\\
\widehat{K} &=\frac{\tau}n\widehat{g}+\phi^{-2}(\sigma+\psi^N\mathring {\mathcal L}_g W).
\end{aligned}\right.
$$
Putting this parametrization in \eqref{eqConstraints}, we  obtain a new interesting system. Of course, we could let $\psi$ depend on
$\phi$ and possibly on some other parameters in the  system  obtained. 
It can be seen  (compare \cite{BartnikIsenberg2004}  for the GR initial data) that
the  conformal method A consists in choosing $\psi=1$, the conformal method B arises when  $\psi=\phi$, and for
$\psi$ a fixed positive function, we obtain  the conformal thin sandwich method.
But many other choices can be made, and some of them  also
yield conformally covariants  systems (see \cite{Delay:paramCC}).

\begin{remarks}
$\mbox{ }$\\
$\bullet$ It is probable that, depending on the gravity theory studied,  an adapted choice of $\psi$
will be judicious.\\
$\bullet$ For an arbitrary  function $f$,  we could use in the same way the York decomposition  related to $\widetilde g$
of the tensor $f( \widehat{K}-\frac{\tau}n\widehat{g})$. This will produce a system where
$f$ and $\psi$ can be choosen freely (possibly depending on some other parameters and/or variables).\\
$\bullet$ In \cite{ToriiShinkaiEGBgravity}, a discussion is made about  a choice of the different powers  of the conformal factor that can be used to define the conformal parametrization.
With  our  parametrizations, only the $\hat A_{ij}$ there is changed.  But we can see that in  that case the natural choice to make, with the notations of
\cite{ToriiShinkaiEGBgravity} is $\tau=0$, $m=2/(N-2)$ and $l=-2$ there.

\end{remarks}

\section{ Appendix }
\subsection{Conformal covariance}\label{cc}
Let us consider three products of  tensor bundles over $M$,
$$
E=E_1\times...\times E_k,\; \;F=F_1\times...\times F_l,\;\;G=G_1\times...\times G_m,
$$
and a differential operator acting on the sections :
$$
P_g:\Gamma(E)\longrightarrow \Gamma(F),
$$
with coefficients  determined 
by $g=(g_1,...,g_m)\in G$. We will say that $P_g$  is conformally covariant
if there exist $a=(a_1,..,a_k)\in\R^k$, $b=(b_1,..,b_l)\in\R^l$ and
$c=(c_1,..,c_m)\in\R^m$ such that for each smooth  section $e$ of $E$, and every smooth function  $\psi$ on
 $M$, we have
$$
\psi^b \odot P_{\psi^c \odot g}(\psi^a \odot e)=P_g(e),
$$
where  
$$
\psi^a\odot e=(\psi^{a_1} e_1,...,\psi^{a_k} e_k).
$$
A differential system will be said conformally covariant if it can be written  in the form
 $P_g(e)=f$,
for a conformally covariant operator $P_g$.

\subsection{The equations}\label{sectionequation}
We specify here the equations related to the Hamiltonian constraint $\rho(g,K)=0$ and the momentum constraint $J(g,K)=0$, in some modified gravity theories.\\

The constraints equations for the general Lovelock gravity can be found in 
\cite{YCBLivre2009} pages 692-693 (a reedition of \cite{YCBGaussBonnet1989} pages 56-57).
We do not reproduce them here but choose to only details the Einstein-Gauss-Bonnet
particular case.\\
$\mbox{  }$\\
{\bf Einstein-Gauss-Bonnet constraint:}
The following expressions can be found in \cite{ToriiShinkaiEGBgravity} for instance:
$$
\rho(g,K)=M+\alpha_{GB}(M^2-4M_{ij}M^{ij}+M_{ijkl}M^{ijkl}),
$$
$$
\frac{-1}{2}J_i(g,K)=N_i+2\alpha_{GB}(MN_i-2M_i^jN_j+2M^{kl}N_{ikl}-M_i\;^{jkl}N_{klj}),
$$
where
$$
M_{ijkl}=R_{ijkl}(g)+K_{ik}K_{jl}-K_{il}K_{jk},
$$
$$
M_{ij}=R_{ij}(g)+\Tr_gKK_{ij}-K_{il}K^l_j,
$$
$$
M=R(g)+(\Tr_gK)^2-K_{ij}K^{ij},
$$
$$
N_{ijk}=\nabla_iK_{jk}-\nabla_jK_{ik}
$$
$$
N_i=\nabla_jK^j_i-\nabla_i\Tr_gK,
$$
and $\alpha_{GB}$ is a coupling constant, equal to zero in the Einstein theory. 

$\mbox{  }$\\
{\bf f(R) gravity constraint:}
The following equations can be found in \cite{LeflochMafR}~:
$$
\begin{aligned}
\rho(g,K)&=R(g)-K_{ij}K^{ij}+(\Tr_g K)^2\\
&-f'(\mathcal R)^{-1}
[2\Delta_g f'(\mathcal R)+2\Tr_gKf''(\mathcal R)\dot{\mathcal R}-f(\mathcal R)-\mathcal Rf'(\mathcal R)],
\end{aligned}
$$
$$
\frac12J_j(g,K)=-\nabla^iK_{ij}+\nabla_j\Tr_gK-f'(\mathcal R)^{-1}\nabla_j(f''(\mathcal R)\dot{\mathcal R})-K^i_j\nabla_i(\ln(f'(\mathcal R))),
$$
where $\Delta=\nabla^i\nabla_i$, and $\mathcal R$, $\dot{\mathcal R}$ are respectively the  initial data\footnote{to be chosen}
of the scalar curvature and the time derivative  of the scalar curvature of the space time.
More precisely $\dot{\mathcal R}$ will be equal to $\mathcal L_n\mathcal R$, where
$n$ will be the futur unit normal to the initial data of the  space-time when evolved.

\def\polhk#1{\setbox0=\hbox{#1}{\ooalign{\hidewidth
  \lower1.5ex\hbox{`}\hidewidth\crcr\unhbox0}}}
  \def\polhk#1{\setbox0=\hbox{#1}{\ooalign{\hidewidth
  \lower1.5ex\hbox{`}\hidewidth\crcr\unhbox0}}} \def\cprime{$'$}
  \def\cprime{$'$} \def\cprime{$'$} \def\cprime{$'$}
\providecommand{\bysame}{\leavevmode\hbox to3em{\hrulefill}\thinspace}
\providecommand{\MR}{\relax\ifhmode\unskip\space\fi MR }
\providecommand{\MRhref}[2]{%
  \href{http://www.ams.org/mathscinet-getitem?mr=#1}{#2}
}
\providecommand{\href}[2]{#2}

\end{document}